\documentclass[a4paper,12pt]{article}
\usepackage{amsmath}
\usepackage{amsthm}
\usepackage{amssymb}
\usepackage{amscd}
 \makeatletter

 \@addtoreset{equation}{section}
  \makeatother

\theoremstyle{plain}
\newtheorem{thm}{Theorem}[section]
\newtheorem{prop}[thm]{Proposition}
\newtheorem{cor}[thm]{Corollary}
\newtheorem{lem}[thm]{Lemma}
\theoremstyle{definition}
\newtheorem{defn}{Definition}[section]

\theoremstyle{remark}
\newtheorem{rem}{Remark}[section]

\newenvironment{namelist}[1]{%
\begin{list}{}
{
\settowidth{\labelwidth}{#1}
\setlength{\leftmargin}{1.1\labelwidth}}
}{%
\end{list}}

\begin{document}
\title{  Four-dimensional Wess-Zumino-Witten actions }
\author{ Tosiaki Kori\thanks
{Research supported by
Promotion for Sciences of the ministry of
education in Japan ( no. 13640224 )
and Waseda University Grant for special research project ( no. 2000A-129)}\\
Department of Mathematics\\
School of Science and Engineering\\
 Waseda University \\
3-4-1 Okubo, Shinjuku-ku
 Tokyo, Japan.\\e-mail: kori@waseda.jp
}
\date{ }
\maketitle
\begin{abstract} 
We shall give an axiomatic construction of  Wess-Zumino-Witten actions 
valued in  \(G=SU(N)\), \(N\geq 3\).    It is realized as a
functor \({WZ}\) from the category of conformally flat four-dimensional
manifolds to the category of line bundles with connection that
satisfies, besides the axioms of a topological field
theory, the axioms which abstract the characteristics of
Wess-Zumino-Witten actions.    To each conformally flat
four-dimensional manifold 
\(\Sigma\) with boundary
\(\Gamma=\partial\Sigma\), a line bundle
\(L=WZ(\Gamma)\) with connection over the space \(\Gamma G\) of 
mappings from \(\Gamma\) to 
\(G\) is associated.   The Wess-Zumino-Witten action is 
 a non-vanishing horizontal section \(WZ(\Sigma)\) 
of the pullback bundle
\(r^{\ast}L\)  over \(\Sigma G\) by the boundary restriction \(r:\Sigma
G\longrightarrow \Gamma G\).  
   \(WZ(\Sigma)\) is required to satisfy  
 a generalized Polyakov-Wiegmann formula with
respect to the pointwise multiplication of the fields \(\Sigma G\).   
  Associated to
the WZW-action there is a geometric description of the extension of the Lie group 
 \(\Omega^3G\) due to  J. Mickelsson.   In fact we have two abelian 
extensions of \(\Omega^3G\)  that are in duality.
\end{abstract}


{\it Subj. Class.:} Global analysis, Quantum field theory

Conformally flat four-manifolds.

\medskip

\section{Introduction}
In this paper we shall give an axiomatic construction of the
 Wess-Zumino-Witten action.   Axiomatic approaches to field
theories were introduced by 
 G. Segal in two-dimensional conformal
field theory (CFT), and by M. F. Atiyah in topological
 field theory,
 [ 1, 16 ].   The axioms  abstract the functorial structure that
the path integral would create if it existed as a mathematical object.   
 Thus a CFT is defined as a Hilbert space representation of the
operation of disjoint union and contraction on a category of manifolds 
 with parametrized boundaries.   The functional
integral formalism was also explored by Gawedzki [ 7 ] to
explain the WZW conformal field theory.    M. A. Singer [ 18 ]
proposed a four-dimensional CFT in the language of Penrose's
twistor space, where Riemann surfaces of
two-dimensional CFT were replaced by conformally flat
four-dimensional manifolds.   
  
  In a
four-dimensional Wess- Zumino-Witten model the space of field
configurations is the space of all maps from closed four-dimensional
manifolds with or without boundary into a compact Lie group. 
We know from the discussions in [ 18, 21 ] that the geometric setting for CFT is most naturally given by the category of conformally flat
manifolds.   So we adopt this category of manifolds also for our
WZW model.    Let \(\Sigma\) be a conformally flat four-dimemsional
manifold with boundary
\(\Gamma=\partial\Sigma\) which may be the empty set.   Let
\(G=SU(N)\) with \(N\geq 3\).   
 The amplitude of the WZW model is given formally by the functional
integration 
 over fields 
\(f\in\Sigma G=Map(\Sigma,G)\) with the boundary restriction equal to
the prescribed \(g\in\Gamma G=Map(\Gamma,G)\):
\begin{equation}
A_{\Sigma}(g)=\int_{f\in\Sigma G;\,
f\vert\Gamma =g}\,\exp\{2\pi iS_{\Sigma}(f)\}\,{\cal D}f,
\end{equation}
where 
\(S_{\Sigma}(f)\) is 
 defined  by;
\begin{equation}
S_{\Sigma}(f)=-\frac{ik}{12\pi^2}\int_{\Sigma}
tr(df^{-1}\wedge\ast df)+C_{\Sigma}(f)\nonumber. 
\end{equation}
Since we deal with contributions that are topological in nature we omit the first 
term ( kinetic term ).   The exponential of the second
term 
\begin{equation}
WZ(\Sigma)(f)=
\exp\{2\pi i C_{\Sigma}(f)\} \end{equation}
is called the Wess-Zumino-Witten action.    ( In [ 7, 8 ] it is
called an amplitude or a probability amplitude.   In [ 3 ] it is
called the Wess-Zumino-Witten action.)    When
\(\Sigma\) has no boundary
\(C_{\Sigma}(f)\) is defined by
\begin{equation}
C_{\Sigma}(f)=\frac{i}{240\pi^3}\int_{B^5} tr(d\tilde{f}\cdot
\tilde{f}^{-1})^5 ,\end{equation}
 where \(\tilde f\) is an extension of \(f\) to a 5-dimensional manifold \(B^5\) with boundary \(\partial
B^5=\Sigma\).   
 Since \(\Sigma\) is a compact conformally flat manifold it is the boundary of a five-dimensional manifold  \(B^5\).   But it is not clear that we can take such a smooth extension of \(f\) over \(B^5\).    If \(\Sigma\) is simply connected it is conformally equivalent to a four-dimensional sphere, and then, since \(\pi_4(G)=1\), there exists a smooth extension of \(f\) to the five-dimensional disc \(D^5\) and \(C_{S^4}(f)\) is defined up to  
\({\bf Z}\), that is, \(\exp\{2\pi iC_{S^4}(f)\}\) is well defined.   The problem arises as to how to define the action \(WZ(\Sigma)(f)\) for general \(\Sigma\) without boundary.     
On the other hand in ( 0.1 ) we are dealing with a four-manifold with
boundary, so we must also give the definition of the action 
 \(WZ(\Sigma)(f)\) 
 for \(\Sigma\) with non-empty boundary.  The above discussions lead to the following conclusion:
{\it A four-dimensional Wess-Zumino-Witten ( WZW ) model means to assign a proper definition of the action 
\(WZ(\Sigma)(f)\) to every compact conformally flat 
four-manifold \(\Sigma\) with or without
boundary.}   

We shall construct the actions \(WZ(\Sigma)\) as the
 objects that satisfy several  
axioms.  
  Our WZW actions are associated to four-dimensional manifolds
 with boundary
and respect the functorial properties of various  
operations on the basic 
manifolds.   Hence we impose on \(WZ(\Sigma)\) several axioms
that are similar to those of topological field theories.  
Axioms of topological field theories were introduced by M. F.
Atiyah in [ 1 ].   They apply to a functor from the category of topological
spaces to the category of vector spaces.    K. Gawedzki
explored in the same spirit the axioms which characterize the
amplitudes of two-dimensional WZW theory, [ 7 ].   Since our
objects are not the amplitudes but the actions of the field,
we describe our four-dimensional WZW theory as a functor
\(WZ\) from the category of four-manifolds with boundary to
the category of {\it complex line bundles}.   This functor is
required to satisfy the involutory axiom, the multiplicativity
axiom and the associativity axiom that represent respectively
the orientation reversal and the operations of disjoint union and
contraction of the basic manifolds.    Next we shall introduce
two axioms that are characteristic of WZW models.   We know that the
action functional in field theory has topological effects,
that is, it gives rise to the holonomy of a connection.   So
we require as our next axiom that the action
\(WZ(\Sigma)\) gives rise to a four-dimensional analogue of
parallel transport associated to a connection of the complex
line bundle.   Higher-dimensional parallel transports as well
as holonomies were discussed by Y. Terashima, [ 19 ],
following the idea of Gawedzki in [ 8 ] that relates
isomorphism classes of line bundles with connection and
the \(U(1)\)-holonomy coming from WZW action.   The
fundamental property of the WZW action is its behavior under the 
pointwise
multiplication of fields.   It is expressed by the
Polyakov-Wiegmann formula, [ 13 ], and its 
generalization to four-dimensional sphere was given by J. Mickelsson,[ 11 ].   As our last
axiom we demand that
\(WZ(\Sigma)\) satisfies the generalized Polyakov-Wiegmann foumula
 over \(\Sigma G\).   
     More precisely the
WZW actions can be stated as follows.    A four-dimensional
WZW model means a functor \(WZ\) that assigns to each
manifold
\(\Sigma\), and its boundary
\(\Gamma=\partial\Sigma\),
 a line bundle \(L=WZ(\Gamma)\) over the space 
of maps \(\Gamma G\), and a non-vanishing section \(WZ(\Sigma)\) 
over \(\Sigma G\)
of the pullback line bundle \(r^{\ast}L\) by the boundary restriction map \(
r:\,\Sigma G\longrightarrow \Gamma G\).   The functor \(WZ\)
satisfies the axioms of topological field theories.  
We demand that each line bundle \(WZ(\Gamma)\) has
a connection and that  
\(WZ(\Sigma)\) is parallel with respect to the induced
connection on
\(r^{\ast}L\).    We impose moreover that on \(r^{\ast}L\)
there is defined a product which is equivariant with respect
to the product on \(\Sigma G\) through the  Polyakov-Wiegmann
formula:
\begin{equation}
WZ(\Sigma)(fg)=WZ(\Sigma)(f)\ast
WZ(\Sigma)(g)\qquad\mbox{for
$f,g\in\Sigma G$.}\end{equation}   
We shall see that \(WZ(\Sigma)\) is a positive integer for a compact \(\Sigma\).

Here is a brief summary of each section.   In section 1,  we explain following [ 18 ] that the category of
conformally flat manifolds fits most
naturally the construction of axiomatic CFT and our
WZW model.   In 1.2 we introduce the axioms of our WZW
model.    Gawedzki in [ 7 ] gave two line bundles in duality 
over the loop space \(LG\) that correspond to the 2-cocycles
obtained by transgressing the  3-curvature on \(G\).   In the same
spirit we shall give in section 2 two line bundles \(WZ(S^3)\) 
and \(WZ((S^3)')\) in duality over
\(\Omega^3_0G\) that correspond to the 2-cocycles obtained by
transgressing the 5-form over \(G\).    Here \(\Omega^3_0G\) is
the space of smooth maps from \(S^3\) to \(G\) that have degree
\(0\).   In fact 
 we have a   two-form on \(\Omega^3_0G\); 
   \begin{equation}\beta=
\frac{i}{240\pi^3}\int_{S^3}\,tr(df\cdot f^{-1})^5 ,\end{equation}
which generates the integral cohomology
class \(H^2(\Omega_0^3G,{\bf Z})\).   Hence it defines a line bundle with
connection on \(\Omega_0^3G\), with the curvature \(\beta\).  
This is \(WZ(S^3)\).   
 Let \(DG\) be the space of maps
from a hemisphere \(D\) to \(G\) and let \(D'G\) be the space of maps for
 the other hemisphere.   We shall give a non-vanishing section \(WZ(D)\) of 
the pullback line bundle of \(WZ(S^3)\) by the boundary restriction map
\(r: DG\longrightarrow \Omega^3_0G\).   Intuitivvely \(WZ(D)(f)\)
is the holonomy associated to the curvature \(\beta\) over the
four-dimensional path \(f\in DG\).   Similarly we have a
non-vanishing section
\(WZ(D')\) of the pullback line bundle of \(WZ((S^3)')\) by \(r':
D'G\longrightarrow
\Omega^3_0G\).   The connections on \(WZ(S^3)\) and \(WZ((S^3)^{\prime})\) are 
given in 2.8,
with respect to which \(WZ(D)\) and \(WZ(D^{\prime})\) are  parallel respectively.   
  In section 3 we construct the functor \(WZ\).    
The line bundle \(WZ(\Gamma)\) is defined as the tensor product
of \(WZ(\Gamma_i)\) for each boundary component
\(\Gamma_i\) parametrized by
\(S^3\), while each \(WZ(\Gamma_i)\) is defined as the pullback of 
\(WZ(S^3)\) or \(WZ((S^3)')\) by the map
\(\Gamma_iG\longrightarrow S^3G\) coming from 
the parametrization.
The non-vanishing section \(WZ(\Sigma)\) of \(r^{\ast}WZ(\Gamma)\) 
is defined from the non-vanishing sections 
\(WZ(D)\) and \(WZ(D')\) by cutting and pasting
 methods and by using the dual relations, i.e. the associativity
axiom.     The connection on \(WZ(\Gamma)\) is induced from those on
\(WZ(S^3)\) and \(WZ((S^3)')\) by a standard procedure.   
   \(WZ\) satisfies the axioms that abstract 
the functorial structure of the WZW actions.   
In particular we have the Polyakov-Wiegmann formula generalized to \(\Sigma G\) for any conformally flat four-manifold \(\Sigma\).   
In section 4
we shall  discuss extensions of the Lie group 
\(\Omega^3_0G\)  .   It is a well known observation
that the two-dimensional WZW action gives a geometric
description of central extensions of the loop group, [ 2, 7 ].  
The 
\(U(1)-\)principal bundle
over \(\Omega^3_0G\)  associated to the 
line bundle 
\(WZ(S^3)\) however does not have any group structure
.   Instead J. Mickelsson in [ 11 ] gave an extension of \(\Omega^3_0 G\) 
by the abelian group \(Map({\cal A}_3,U(1)) \), where \({\cal A}_3\) is the space
of connections on \(S^3\).
 We shall explain two 
extensions of Mickelsson's type that are dual to each other.

\section{  Axioms for a 4-dimensional WZW model}

\subsection{}  
   The basic components of four-dimensional CFT are some well
behaved class of four-dimensional manifolds \(M\)
 with parametrized boundaries, together with the natural
operations of disjoint union
\[(M_1,M_2)\longrightarrow M_1\cup M_2,\]
and contraction 
\[M\longrightarrow \tilde M,\]
where \( \tilde M\) is obtained from \(M\) using the parametrization to
 attach a pair of boundary three-spheres to  
 each other.    A four-dimensional CFT is
then defined as
 a Hilbert space representation of the operation of disjoint union and
contraction on these basic components.    
Now we know that the geometric setting for this CFT is most
naturally given by the conformal equivalence classes of 
conformally flat four-dimensional manifolds.   This fact
was explained by M. A. Singer [ 18 ] , R. Zucchini [ 21 ]
and Mickelsson-Scott [ 12 ].   

Here we shall see following [ 18 ] the fact that the class of
compact conformally flat four-dimensional manifolds with boundary
is closed under the operation of sewing manifolds together across a
boundary component.   For any conformally flat \(M\)
the developing map \(M\longrightarrow S^4\) is
a well defined conformal local diffeomorphism.   A
closed 3-manifold \(N\subset M\) is called a
{\it round} \(S^3\) in \(M\) if it goes over
diffeomorphically to a round \(S^3\) in \(S^4\)
under development.   This is well defined because the
developing map is unique up to composition with
conformal transformations.   For standard \(M\), the boundary 
\(\partial M\) consists of a disjoint union of round
\(S^3\)Õs, [ 15 ].   For each boundary component
\(B\) one can find a neighborhood of \(B\) in \(M\)
and a conformal diffeomorphism of this neighborhood
onto a neighborhood of the equator in the northern
hemisphere of \(S^4\).   If we have two boundary
components \(B\) and \(\tilde B\) of \(M\) and an
orientation reversing conformal diffeomorphism
\(\psi: B\longrightarrow \tilde B\), then \(B\) and
\(\tilde B\) can be attached using \(\psi\) and 
the resulting manifold will have a unique
conformally flat structure compatible with the
original one on \(M\).

\subsection{}  
Now we give the precise definition of a 
four-dimensional WZW model.     

Let \({\cal M}_4\) be the conformal equivalence classes of all compact
conformally flat four-dimensional manifolds 
\(M\) with boundary \(\partial M=\bigcup_{i\in
I}\Gamma_i\) such that each oriented component \(\Gamma_i\) is a
round \(S^3\), and is endowed with a parametrization
\(p_i:\,S^3\longrightarrow
\Gamma_i\).   We  
 distinguish positive  and negative parametrizations \(p_i
: S^3\longrightarrow \Gamma_i\,,
\,i\in I_{\pm}\), depending on whether $p_i$ respects 
the orientation of  $\Gamma_i$ or not.    

Let \({\cal M}\) be the category whose objects are 
three-dimensional manifolds \(\Gamma\) which are disjoint unions of round
\(S^3\)Õs.   A morphism between three-dimensional manifolds
\(\Gamma_1\) and \(\Gamma_2\) is an oriented cobordism given by
\(\Sigma\in  {\cal M}_4\) with boundary   
\(\partial\Sigma=\Gamma_2\bigcup(\Gamma_1^{\prime})\), where
the upper prime indicates the opposite orientation.  

Let \({\cal L}\) be the category of complex line bundles.

Let \(G=SU(N)\),\(\,N\geq 3\).   In the following the set of
smooth mappings from a manifold \(M\) to \(G\) that are
based at some point \(p_0\in M\) is  denoted by
\(MG=Map(M, G)\).   \(MG\) becomes a group under product of
mappings.     For a \(\Sigma\in{\cal M}_4\) with boundary
\(\Gamma=\partial\Sigma \), 
\(r\) denotes the restriction map
\begin{equation}
r:\,\Sigma G\longrightarrow \Gamma G,\qquad r(f)=f\vert\Gamma.
\end{equation}

 A four-dimensional WZW model means a functor \(WZ\) from the category
\({\cal M}\) to the category \({\cal L}\) which
assigns:
\begin{namelist}{axiom 1}
\item[{\bf WZ1},]
 to each manifold $\Gamma\in{\cal M}$, a complex line bundle
$WZ(\Gamma)$  over the space $\Gamma G$,
\item[{\bf WZ2},]
to each $\Sigma\in {\cal M}_4$, with  
\(\partial\Sigma=\Gamma$, a non-vanishing section 
\( WZ(\Sigma)\) of the pullback line bundle \(r^{\ast} WZ(\Gamma)\) .
\end{namelist}

Recall that the pullback bundle is by definition
\begin{equation}
r^{\ast}WZ(\Gamma)=\left\{(f,u)\in \Sigma G\times
WZ(\Gamma);\quad \pi u=r(f),\right\} \end{equation}
and the section \(WZ(\Sigma)\) is given at \(f\in\Sigma
G\) by
\[
WZ(\Sigma)(f)=(f,u)\quad\mbox{with } u\in
\pi^{-1}(r(f))=WZ(\Gamma)_{r(f)}.\] 

\(WZ\) being a functor  from 
\({\cal M}\) to \({\cal L}\), a conformal
diffeomorphism \(\alpha:\Gamma_1\longrightarrow \Gamma_2\)
induces an isomorphism \(WZ(\alpha):
WZ(\Gamma_1)\longrightarrow WZ(\Gamma_2)\) such that 
\(WZ(\beta\alpha)=WZ(\beta)WZ(\alpha)\) for
\(\beta:\Gamma_2\longrightarrow \Gamma_3\).   Also if \(\alpha\)
extends to a conformal diffeomorphism
\(\Sigma_1\longrightarrow\Sigma_2\), 
with \(\partial\Sigma_i=\Gamma_i\), \(i=1,2\), then \(WZ(\alpha)\)
takes \(WZ(\Sigma_1)\) to \(WZ(\Sigma_1)\).

The functor $WZ$ satisfies the following axioms.   {\bf A1},{\bf A2} and
{\bf A3} represent in the category of line bundles the orientation reversal and 
the operation of disjoint union and contraction.  
 These axioms are stated in the
same manner as in topological field theories, [ 1 ].   Axioms {\bf A4} and {\bf
A5} are characteristic of the WZW model.

\begin{namelist}{ABC}
\item[{\bf A1}]( Involution ): 
 \begin{equation}WZ(\Gamma^{\prime})=WZ(\Gamma)^{\ast}
\end{equation}
 where $\ast$ 
indicates the dual line bundle.
\item[{\bf A2}]( Multiplicativity ):   
\begin{equation}WZ(\Gamma_1\cup\Gamma_2)=WZ(\Gamma_1)\otimes
WZ(\Gamma_2).\end{equation}
\hskip 0.5cm
\item[{\bf A3}]( Associativity ):

For a composite cobordism 
 \(\Sigma=\Sigma_1\cup_{\Gamma_3}\Sigma_2\) such that 
 \(\partial\Sigma_1=\Gamma_1\cup\Gamma_3\) and 
 \(\partial\Sigma_2=\Gamma_2\cup\Gamma_3^{\prime}\), we have
\begin{equation}WZ(\Sigma)(f)=<WZ(\Sigma_1)(f_1),
WZ(\Sigma_2)(f_2)>,\end{equation}
for any \(f\in\Sigma G\), \(f_i=f\vert\Sigma_i\),
i=1,2,  where
\(<\,,\,>\) denotes the natural pairing
\begin{equation}WZ(\Gamma_1)\otimes WZ(\Gamma_3)\otimes
WZ(\Gamma_3^{\prime})
\otimes WZ(\Gamma_2)\longrightarrow WZ(\Gamma_1)\otimes
WZ(\Gamma_2).\end{equation}

\end{namelist}
More precisely, let 
\(WZ(\Sigma_1)(f_1)=( f_1, u_1\otimes v)
\) and \( WZ(\Sigma_2)(f_2)=(f_2, u_2\otimes
v')\)  with \(u_i\in
WZ(\Gamma_i)\) for \(i=1,2\), and \(v\in WZ(\Gamma_3))\)
,  \(v'\in WZ(\Gamma'_3)\).    From the definition  
\(u_i\in\pi^{-1}(f_i\vert\Gamma_i)\),
\(v\in\pi^{-1}(f_1\vert\Gamma_3)\) and 
\(v'\in\pi^{-1}(f_2\vert\Gamma'_3)\).   
On the other hand, let 
\(WZ(\Sigma)(f)=(f,w_1\otimes w_2)\in
WZ(\Gamma_1)\otimes WZ(\Gamma_2)\) with $w_i\in
\pi^{-1}(f\vert\Gamma_i)$, 
for $i=1,2$.   
Then axiom {\bf A3} says that 
\(
w_1\otimes w_2=<v',v>u_1\otimes u_2\).
\medskip
The multiplicative axiom {\bf A2} asserts that if
\(\partial\Sigma=\Gamma_2\bigcup (\Gamma_1^{\prime})\) ,
 then \(WZ(\Sigma)\) is a section of 
\begin{equation}
r_1^{\ast}WZ(\Gamma'_1)\otimes
r_2^{\ast}WZ(\Gamma_2) =Hom(
r_1^{\ast}WZ(\Gamma_1),r_2^{\ast}WZ(\Gamma_2)).
\end{equation}
Therefore any cobordism \(\Sigma\) between \(\Gamma_1\) and
\(\Gamma_2\) induces a homomorphism of sections of pullback line
bundles
\begin{equation}
WZ(\Sigma):C^{\infty}(\Sigma,r_1^{\ast}WZ(\Gamma_1))\longrightarrow
C^{\infty}(\Sigma,r_2^{\ast}WZ(\Gamma_2)).\end{equation}

\medskip
   
We impose:
\begin{enumerate}
\item
\begin{equation}WZ(\phi)={\rm C}\quad\mbox{for \(\phi\) the empty
3-dimensional manifold,}\end{equation}
\item 
\begin{equation}WZ(S^4)=1\end{equation}
\item
\begin{equation}WZ(\Gamma\times [0,1])=Id.(WZ(\Gamma)\longrightarrow
WZ(\Gamma)).\end{equation}
\end{enumerate}

\begin{cor} If $\Sigma$ has no boundary \(\left(\partial\Sigma=\phi\right)\), 
then 
 $WZ(\Sigma)\in {\rm C}$ .\end{cor}

The following axioms are characteristic of WZW models.

\begin{namelist}{ABC}
\item[{\bf A4}] For each \(\Sigma\in {\cal M}_4\) with 
\(\Gamma=\partial\Sigma\), \(WZ(\Gamma)\) has a
connection, and \(WZ(\Sigma)\) is parallel with respect to the
induced connection on \(r^{\ast}WZ(\Gamma)\) .
\item[{\bf A5}]  ( Generalized Polyakov-Wiegmann formula ):
For each \(\Sigma\in {\cal M}_4\) with \(\Gamma=\partial\Sigma\),
on the pullback line bundle
\(r^{\ast}WZ(\Gamma)\) is defined a product \(\ast\) with respect to which we
have  
\begin{equation}
WZ(\Sigma)(fg)=WZ(\Sigma)(f)\ast WZ(\Sigma)(g)\qquad\mbox{for any $f,g\in\Sigma
G$}.
\end{equation}
\end{namelist}
The well known Polyakov-Wiegmann formula extended by 
J. Mickelsson [ 11 ] is
concerned with the case of the four-dimensional sphere, 
\(\Sigma=S^4\).

\medskip
 From now on we shall construct the functor $WZ$ step by
step. In section 2.5 we shall construct two 
line bundles over \(S^3G\), which are \(WZ(S^3)\) and 
\(WZ((S^3)^{\prime})\).   In section 3 we give the functor
 \(WZ\) of WZW actions step by step starting from  \(WZ(S^3)\) and 
\(WZ((S^3)^{\prime})\).

\section{
Line bundles on $\Omega^3G$}

\subsection{}  
In the following we denote by $\Omega^3G$, instead of \(S^3G\), 
 the set of smooth mappings $f$ from a $S^3$ to $G=SU(N)$ that
are based, i.e., $f(p_o)=1$, at some point $p_o\in S^3$.   
It is known that $\Omega^3G$ is not connected and is
divided into  denumerable sectors labelled by the
soliton number ( the mapping degree ).    Here we
follow the explanation due to I. M. Singer of these facts [ 17
], see also [ 4, 10 ].   Let the  evaluation map,
$ev:\,S^3\times
\Omega^3G\longrightarrow G$, be  defined by
$ev(m,\varphi)=\varphi(m), \,m\in S^3,\,\,\varphi\in \Omega^3G$ 
.   The Maurer-Cartan form $g^{-1}dg$ on G gives the
identification  of the tangent space $T_eG$ at $e\in G$ and
$Lie\,G=su(N)$.    The primitive generators of the cohomology 
$H^{\ast}(G,{\rm R})$ are given by 
\begin{equation}\omega_3=-\frac 1{4\pi^2} tr (g^{-1}dg
)^3,\quad   
 \omega_5=\frac {-i}{2\pi^2} tr (g^{-1}dg )^5,
\quad \cdots.\end{equation}  
Integration on $S^3$ of the pull back of $\omega_{2k-1}$ 
by the evaluation map $ev$ gives us the following $2(k-2)$ form 
on $\Omega^3G$; \begin{equation}
\nu_{2k-1}=(\frac 1{2\pi i})^k\frac{((k-1)!)^2}{(2k-1)!}\int_{S^3}
 tr(\, d\varphi\, \varphi^{-1}\,)^{2k-1},\quad 3\leq 2k-1\leq 2N-1.
\end{equation} 

In particular  $\nu_3$ is the mapping degree of $\varphi$;
\begin{equation}\deg\,\varphi=\frac{i}{24\pi^2}\int_{S^3}\,
tr(d\varphi\,\varphi^{-1})^3.\end{equation} 

\begin{prop}
 \begin{enumerate}
 \item $$
S^3
 Lie\,G\,
 \stackrel{\,\exp\,}{\longrightarrow} \,\Omega^3G\,
\stackrel{\,\deg\,}{\longrightarrow} \,{\rm Z}\longrightarrow
0$$ is exact.
 \item  $$\deg\varphi_1\cdot\varphi_2=\deg\varphi_1+\deg\varphi_2.$$
 \end{enumerate}
\end{prop} 
See [ 4, 10 ].

\subsection{}  

Let \(P_G\) be a \(G\)-principal bundle over \(S^4\).   
Let ${\cal A}$ be the space of connections on \(P_G\), that are 
$Lie\,G$-valued one-forms on \(P_G\).    Let ${\cal G}=S^4G$ be the group of based gauge transformations.    The action of ${\cal G}$ on ${\cal
A}$ is given by
$A_g=g^{-1}Ag+g^{-1}dg$ for $A\in {\cal A}$ and $g\in
{\cal G}$.   $F=F(A)=dA+A^2$ denotes the curvature two-form
of $A$.   

The Chern-Simons form on \(P_G\) is 
\begin{equation}
\omega^0_5(A)=tr\,(AF^2-\frac12
A^3F+\frac1{10}A^5).
\end{equation}
 We have then
$tr(F^3)=d\omega^0_5(A)$.

From [ 22 ] we know the relation
\begin{equation}
\omega^0_5(A_g)-\omega^0_5(A)=d\alpha_4(A;g)+\frac1{10}
tr(dg\cdot g^{-1})^5,\nonumber\end{equation}
with
\begin{equation}
\alpha_4(A;g)=tr[\,-\frac12V(AF+FA-A^3)+
\frac14(VA)^2+\frac12V^3A\,],
\end{equation}
where $V=dg\cdot g^{-1}$.

Let \(D^5\) be
a five dimensional disc with boundary \(\partial
D^5=S^4\). Integration over \(D^5\) 
gives us  the {\it gauge anomaly }:
\begin{eqnarray}
\Gamma(A,g)&=&\frac{i}{48\pi^3}\int_{S^4}
tr\lbrack -V(AF+FA-A^3)+\frac12(VA)^2+V^3A\rbrack+C_5(g),
\nonumber\\
C_5(g)&=&\frac{i}{240\pi^3}\int_{D^5} tr(dg\cdot
g^{-1})^5 ,
\end{eqnarray}
here $g\in S^4 G$ is extended to 
$D^5G$, in fact, we have such an extension by virtue of    
$\pi_4(G)=1$.  $C_5(g)$ may depend on the extension 
 but it can be shown that the difference 
of two extensions is an integer, and
$\exp(2\pi i C_5(g)\,)$ is independent of the
extension.

We put, for $f,\,g\in S^4 G$,
\begin{eqnarray}
\lefteqn{\gamma(f,g)=\frac{i}{24\pi^3}
\int_{S^4}\alpha_4(f^{-1}df,\,g)}\nonumber\\
&&=\frac{i}{48\pi^3}\int_{S^4}
tr\lbrack (dgg^{-1})(f^{-1}df)^3+\frac12(dgg^{-1}f^{-1}df)^2+
\nonumber\\
\nonumber\\
&&\mbox{  }\qquad\mbox{ }+(dgg^{-1})^3(f^{-1}df)\rbrack.\\
\nonumber\\
\mbox{and}&&\nonumber\\
&&\omega(f,g)=
\Gamma( f^{-1}df,g)=\gamma(f,g)+C_5(g).
 \end{eqnarray}

\medskip
\begin{rem}
Here we shall look at Mickelson's 2-cocycle for
his abelian extension of \(\Omega^3G\), [ 11 ].   
 The cochain  $\alpha_4$ in ( 2-5 ) is a one-cochain on the group
 $S^4G$, valued in 
$Map({\cal A}_4, {\rm R})$.   The coboundary
$\delta\alpha_4$ is  given by 
\begin{eqnarray*}
\delta\alpha_4(A:\, g_1,g_2) &=&d\beta+\alpha_4(g_1^{-1}dg_1;\,
g_2)\\
\beta(A; \,g_1,g_2) &=&
-tr\lbrack \frac12(dg_2g_2^{-1})(g_1^{-1}dg_1)
(g_1^{-1}Ag_1)-\frac12(dg_2g_2^{-1})(g_1^{-1}Ag_1)
(g_1^{-1}dg_1)\rbrack .
\end{eqnarray*}
Mickelson's 2-cocycle \(\gamma_{\Delta}(A;f,g)\) is defined as the
integration of 
\(\delta\alpha_4(A;g_1,g_2)\) over any region \({\Delta}\subset
S^4\):
\begin{equation}\gamma_{\Delta}(A;f,g)=
\frac{i}{24\pi^3}\int_{\Delta}\delta\alpha_4(A;f,g).
\end{equation}
But for \({\Delta}=S^4\) it is independent of \(A\) and 
\begin{eqnarray}\gamma_{S^4}(A;f,g)&=&\int_{S^4}\delta\alpha_4(A;f,g)\nonumber\\
&=&\int_{S^4}\alpha_4(f^{-1}df,g)=\gamma(f,g),\end{eqnarray}
for  \(f,g\in S^4G\).  Hence, instead of \(\gamma_{S^4}(A;f,g)\), we use
more simple 
\(\gamma(f,g)\) for our purpose.    
\end{rem}
\begin{rem}
We have 
\begin{equation}
\gamma(F,G)=\gamma_D(A;F,G)+\gamma_{D^{\prime}}(A;F,G),
\end{equation}
for any \(A\in {\cal A}_4\).   Here \(D\) is an oriented hemisphere of \(S^4\) and \(D^{\prime}\) is the other hemisphere: \(D\cup D^{\prime}=S^4\).
\end{rem}

\begin{lem}[Polyakov-Wiegmann]    For $f,\,g\in
S^4 G$ we have
\begin{equation}
C_5(fg)=C_5(f)+C_5(g)+\gamma(f,g)\quad
\mod{\bf Z}.\end{equation}
\end{lem}
The following formula was proved by Mickelsson in Lemma 4.3.7 of his book
[ 10 ] .   
\begin{equation}
C_5(fg)=C_5(f)+C_5(g)+\gamma_{S^4}(A;f,g)\quad
\mod{\bf Z}.\nonumber\end{equation}
Since \(\gamma_{S^4}(A;f,g)=\gamma(f,g)\) from ( 2.10 ) we have
the proposition.

\subsection{}  
Now we are prepared to define the line bundle
$WZ(\phi)$ over
  $Map(\partial S^4,G)=\phi$, and the section \(WZ(S^4)\) of
the pullback line bundle of \(WZ(\phi)\) by the empty
restriction map \(r:S^4G\longrightarrow \phi\).   

Let $L_{\phi}$ be the quotient of  
$S^4 G\times {\rm C}$ by the equivalence relation;
\begin{equation}(f,c)\,\sim\, (g, c\exp\{2\pi i\omega(f,f^{-1}g)\}\,).
\end{equation}
Then \(L_{\phi}\) is a line bundle
over 
\(Map(\partial S^4,G)=\phi\) with the transition function \(\exp\{2\pi i\omega(f,f^{-1}g)\}\), which we shall define as \(WZ(\phi)\).  Recall that  \(S^4G\) is contractible.   We have
then
\begin{equation}WZ(\phi)\simeq{\rm C}.\end{equation}
The isomorphism is given by 
\[ [f,c]\longrightarrow c\exp\{-2\pi iC_5(f)\}.\]
It is well defined because of the Polyakov-Wiegmann formula.
 Let \(r^{\ast}WZ(\phi)\) be 
the pullback line bundle of \(WZ(\phi)\) by the empty
restriction map \(r:S^4G\longrightarrow \phi\). 
  The section \(WZ(S^4)\) of  \(r^{\ast}WZ(\phi)\) over any \(f\in S^4G\)
is given by \begin{equation} WZ(S^4)(f)=[f,\exp\{2\pi
iC_5(f)\}]\in WZ(\phi).\end{equation}

  By the
isomorphism of (2.14) we can also write
\begin{equation} WZ(S^4)=1\in{\rm C}.\nonumber\end{equation}

We can define the product on the line bundle  \(WZ(\phi)\simeq{\rm C}\) in an obvious way, but we shall look this product more precisely, rather superfluously, for the sake of later sections.    
In \(S^4G\times {\rm C}\) we define the product by
putting;
\begin{equation}
(f,a)\,\ast\,(g,b)=(fg,\, ab\exp\{ 2\pi i
\gamma(f,g)\}\,).\end{equation}
Since the equivalence relation ( 2.13 ) respects the product, it gives a
product on the line bundle \(WZ(\phi)\).   The Polyakov-Wiegmann formula (
2.12 ) is stated as follows.
\begin{equation}
WZ(S^4)(fg)=WZ(S^4)(f)\ast WZ(S^4)(g)\qquad\mbox{ for $f,g\in S^4G$}.
\end{equation}

\subsection{}  

In this paragraph we shall prepare some notations,
definitions and elementary properties that will be used in
the following sections.

Let \(\Omega^3G\) be as before the set of smooth mappings
from
\(S^3\) to \(G=SU(N)\) that are based.
$\Omega^3G$ is not connected but divided into  the
connected components by
$\deg$.   We put 
\begin{equation}\Omega^3_0G=\{
g\in \Omega^3G;\, \deg g=0\}.\end{equation}

The oriented 
4-dimensional disc 
with boundary $S^3$ is denoted by $D$, while that with oposite
orientation 
is denoted by $D^{\prime}$.   The composite cobordism of $D$
and $D^{\prime}$ becomes $S^4$.   We write as before $DG=Map(D,G)$ and 
$D^{\prime} G=Map(D^{\prime}, G)$.   The restriction to $S^3$ of a 
$f\in DG$ has degree 0; $\,f\vert S^3\in \Omega^3_0G$.  

\medskip

For an $a\in \Omega^3_0G$ we denote by $Da$ the set of those
$g\in DG$ that is a smooth extension
of $a$, respectively \(D'a\) is the set of those \(g'\in D'G\)
that is a smooth extension of \(a\).    For
$f\in Da$ and
$g\in Db$ one has $fg\in D(ab)$,   and every element of $D(ab)$ is
of this form.   Similarly for \(D^{\prime}(ab)\).    We
denote by
$g\vee g^{\prime}\in S^4G$ the map obtained by sewing
$g\in DG$ and $g^{\prime}\in D^{\prime}(g\vert S^3)$.

 The upper prime will indicate  that the function expressed by
the letter  is defined on $D^{\prime}$, for example, $1^{\prime}$ is the 
constant function $D^{\prime}\ni x\longrightarrow 1'(x)=e\in G$, while 
 $1$ is the constant function 
\(\, D\ni x\longrightarrow 1(x)=e\in
G\).  We write 
\[
D'f=\{f'\in D'G:\, f'\vert S^3=f\vert S^3\},\qquad
Df'=\{f\in DG:\, f\vert S^3=f'\vert S^3\}.\]

\medskip

Let \(f, g\in DG\) and \(f\vert S^3=g\vert S^3\).   From ( 2.7 )
and ( 2.8 ) we see that
\(\gamma(f\vee f',f^{-1}g\vee 1')\) and  \(\omega(f\vee
f',f^{-1}g\vee 1')\)  are independent of \(f'\in Df\);
\begin{eqnarray}
\gamma(f\vee f',f^{-1}g\vee 1')&=&
\frac{i}{48\pi^3}\int_{D}
tr\lbrack (dgg^{-1})(f^{-1}df)^3+\frac12(dgg^{-1}f^{-1}df)^2+
\nonumber\\
\nonumber\\
&&\mbox{  }\qquad\mbox{ }+(dgg^{-1})^3(f^{-1}df)\rbrack.
\end{eqnarray}
 Similarly,
for \(f', g'\in D'G \) such that \(f'\vert S^3=g'\vert S^3\).
 \(\gamma(g\vee g',1\vee
(g')^{-1}f')\) and  \(\omega(g\vee g',1\vee
(g')^{-1}f')\)  are independent of \(g\in Dg'\). 
   Hence
    
\noindent \(\exp\{2\pi i\omega(f\vee\cdot \,,f^{-1}g\vee 1')\}\) 
and
\(\exp\{2\pi i\omega(\cdot\,\vee f',1\vee
(f')^{-1}g')\}\) 
are constants of \(U(1)\).

\begin{defn}
\begin{enumerate}
\item
We put, for \(f,g\in DG\) such that \(f\vert S^3=g\vert S^3\), 
\begin{equation}
\chi(f,g)=\exp\{2\pi i\omega(f\vee
\cdot \,,f^{-1}g\vee 1')\}.
\end{equation}
\item
We put, for \(f',g'\in D'G\) such that \(f'\vert S^3=g'\vert S^3\), 
\begin{equation}
\chi^{\prime}(f',g')=
\exp\{2\pi i\omega(\,\cdot\vee f',1\vee
(f')^{-1}g')\}.
\end{equation}
\end{enumerate}
\end{defn}

\begin{lem}
\begin{enumerate}
\item
For \(f,\,g\in DG\) such that \(f\vert S^3=g\vert S^3\), we have
\(\chi(f,g)\in U(1)\) and 
\begin{equation}
\exp\{2\pi iC_5(g\vee f')\}=\exp\{2\pi iC_5(f\vee f')\}\chi(f,g)
\qquad
\mbox{ for any $f'\in D'f$}
\end{equation}
\item
For \(f',\,g'\in D'G\)  such that \(f'\vert S^3=g'\vert S^3\),we have
\(\chi'(f',g')\in U(1)\) and 
\begin{equation}
\exp\{2\pi iC_5(f\vee g')\}=\exp\{2\pi iC_5(f\vee f')\}\chi'(f',g'),\qquad
\mbox{ for any $f\in Df'$}
\end{equation}
\end{enumerate}
\end{lem}
The lemma follows from the Polyakov-Wiegmann formula.

\subsection{ }
Now we shall give two line bundles on \(\Omega^3_0G\) that are
 dual to each other.  We shall follow the  arguments due to K. Gawedzki [ 7 ], that were developed to
construct two line bundles in duality over the loop group \(LG\)
and to give the definition of WZW action on a hemisphere.   

We consider the following quotient;
\begin{equation}
 L=D'G\times {\rm C}/\sim' ,
\end{equation}
where \(\sim'\) is the equivalence relation defined by 
\begin{equation}(f',c')\sim' (g',d') \mbox{ if and only if }
\left\{
\begin{array}{rl} f'\vert S^3&=g'\vert S^3\\
d' &=c'\chi'(f',g')
\end{array}
\right. \end{equation}
The equivalence class of \((f',c')\) is denoted by \([f',c']\).   
We define the projection \[\pi: L\longrightarrow \Omega^3_0G\] by 
$\pi([f',c' ])=  f'\vert S^3$.   \(L\) becomes a line bundle on \(\Omega^3_0G\) with the transition function \(\chi'(f',g')\).

More precisely,
let \(a\in \Omega^3_0G\) and take  
\(f'\in D'a\).   
A coordinate neighborhood of \(a\) is given by 
\begin{eqnarray}
U_{f'}&=&\{g'\vert S^3; \quad g'\in V_{f'}\}\nonumber\\
V_{f'}&=&\{g'\in D'G, \,\,g'=\exp X\cdot f'\,;\quad X\in
D'(Lie\,G),\,\Vert X\Vert<\delta\,\}.\nonumber
\end{eqnarray} 
The local trivialization of \( L\) 
is given by the map \(\pi^{-1}(U_{f'})\ni [h',c']\longrightarrow
(h'|S^3,\,c')\);
\begin{equation}
\pi^{-1}(U_{f'})\simeq U_{f'}\times {\rm C}.\nonumber\end{equation}
The transition function \(\chi_{U_{f'},U_{g'}}(b)\) of \( L\) at
 \(b\in U_{f'}\cap U_{g'}\) becomes as follows.   Let \(b\in U_{f'}\cap U_{g'}\).   Let \(h'\in V_{f'}\) and \(k'\in V_{g'}\) be such that \(h'\vert S^3=k'\vert S^3=b\).   For \(\xi=[h',c']=[k',d']\in \pi^{-1}(b)\) we
have  obviously \(d'=\chi'(h',k')c'\).   Hence 
\begin{equation}
\chi_{U_{f'},U_{g'}}(b)=\chi'(h',k')\end{equation}

The line bundle \( L\) is what we wanted to construct and will be denoted by
\(WZ(S^3)\). 

In regard to the involution axiom {\bf A1} which
\(WZ(\cdot)\) is required to satisfy 
we must define another line bundle on \(\Omega^3_0G\) corresponding to   
\(S^3\) with opposite orientation.   
This line bundle \(WZ((S^3)')\) is defined by 
\begin{equation}
WZ((S^3)')=DG\times {\rm C}/\sim\end{equation}
with the equivalence relation
\begin{equation}(f,c)\sim (g,d) \mbox{ if and only if }
\left\{
\begin{array}{rl} f\vert S^3&=g\vert S^3\\
d &=c\chi(f,g)
\end{array}
\right. \end{equation}
 The projection \(\pi:WZ((S^3)')\longrightarrow \Omega^3_0G\) is given by 
 \([f,c]\longrightarrow f|S^3\).   It is a line bundle with the transition function
\(\chi(f,g)\).

 \(WZ(S^3)\) and \(WZ((S^3)')\) are in duality so that the involution 
axiom {\bf A1} is verified for these line bundles.   In fact, the
duality
\[ WZ(S^3)\times WZ((S^3)^{\prime})\longrightarrow
{\rm C}\]
 is defined by 
\begin{equation}<\,[f',\,c'\,],\,[f,\,c\,]>
=cc'\exp\{-2\pi i C_5(f\vee f^{\prime})\},
\end{equation}
where $f\vert S^3=f^{\prime}\vert S^3\in \Omega^3_0G$.   
If we note the evident fact that
 $\gamma(F,1\vee h^{\prime})$ ( resp. $\gamma(F, h\vee 1')$ )
 in ( 2.19 ) 
 is given by an integration over $D^{\prime}$ ( resp. $D$ ), we 
see that the product of transition rules 
( 2.25 ) and ( 2.28 ) imply the transition rule ( 2.13 ) 
of $WZ(\phi)$;
\begin{equation}
\chi(f,g)\chi'(f',g')=\exp\{2\pi i\omega(f\vee
f',f^{-1}g\vee(f')^{-1}g')\},\end{equation}
Hence \begin{equation}
WZ(S^3)\otimes WZ((S^3)')=WZ(\phi).\end{equation}
Composed with ( 2-14 ) this implies the above duality.

 \subsection{ }
 
Let \( r:DG\longrightarrow S^3G\) and \( r':D'G\longrightarrow S^3G\)
 be the restriction maps.  
 
 We put, for \(f\in DG\),
 \begin{equation}
 WZ(D)(f)=\,[f',\exp\{2\pi iC_5(f\vee f')\}\,]\in WZ(S^3)\vert_{r(f)}.
 \end{equation}
Then we see from Lemma 2.3 that \(WZ(D)\) gives a non-vanishing section of
 the pullback line bundle  \(r^{\ast}WZ(S^3)\)
.

In the same way we put, for \(f'\in D'G\),
\begin{equation} WZ(D')(f')=[\,f,\exp\{2\pi
iC_5(f\vee f')\}\,]\in WZ((S^3)')\vert_{r'(f')}. \end{equation}
\(WZ(D')\) defines a non-vanishing section of
\((r')^{\ast}WZ((S^3)^{\prime})\).

\begin{prop}
For \(f\in DG\) and \(f'\in D'G\) such that \(f\vert S^3=f'\vert S^3\).  
\begin{equation}
<\,WZ(D)(f),\,WZ(D')(f')\,>=
WZ(S^4)(f\vee f')
\end{equation}
\end{prop}
In fact both sides are equal to \(\exp\{2\pi iC_5(f\vee f')\}\).
\medskip

 \subsection{ }

The total space of the pullback bundle  \(r^{\ast}WZ(S^3)\) 
is written as
\begin{equation}
r^{\ast}WZ(S^3)=
\left\{(f,\lambda);\, f\in DG,\,\lambda=[f',c']\in WZ(S^3)_{r(f)}\,\right\}.
\nonumber\end{equation}
We define the product in  \(r^{\ast}WZ(S^3)\) by the formula;
\begin{equation}(f,\lambda)\,\ast\,(g,\mu)=(fg,\nu),\end{equation}
where, for \(\lambda=[f',a']\in WZ(S^3)_{r(f)}\) and \(\mu=[g',b']\in WZ(S^3)_{r(g)}\), 
  \(\,\nu=[f'g',c']\in WZ(S^3)_{r(fg)}\,\) is defined by
\begin{equation}
c'=a'b'
\exp\{2\pi i\gamma(f\vee f^{\prime},g\vee g^{\prime})\}.
\end{equation} 
 \(\nu\) does not depend on the representations of \(\lambda\) and \(\mu\),
and the product is well defined.    

We have
\begin{equation}
WZ(D)(fg)=WZ(D)(f)\ast WZ(D)(g)\qquad\mbox{ for $f,g\in DG$.}\end{equation}
In fact, this follows from the definition 
\[ WZ(D)(f)=[f',\exp\{2\pi iC_5(f\vee f')\}]\]
and the Polyakov-Wiegmann formula.  

Similarly we have the product on \((r')^{\ast}WZ((S^3)')\) over \(D'G\).   
It is given by 
\begin{equation}
(f^{\prime},\alpha)\,\ast(g^{\prime},\beta)
=(f^{\prime}g^{\prime},\gamma),\end{equation}
where, for \(\alpha=[f,a]\in WZ((S^3)')_{r'(f')}\) and \(\beta=[g,b]\in WZ((S^3)')_{r'(g')}\),  \(\gamma=[fg,c]\in WZ((S^3)')_{r'(f'g')}\,\) is defined
by
\begin{equation}
c=ab
\exp\{2\pi i\gamma(f\vee f^{\prime},g\vee g^{\prime})\}.
\end{equation}

We have
\begin{equation}
WZ(D')(f'g')=WZ(D')(f')\ast WZ(D')(g')\qquad\mbox{ for $f',g'\in
D'G$.}
\end{equation}

We note that product operations on
\(r^{\ast}WZ(S^3)\) and on 
\((r')^{\ast}WZ((S^3)')\) are compatible with the duality
\begin{equation}r^{\ast}WZ(S^3)\times
(r')^{\ast}WZ((S^3)')\longrightarrow WZ(\phi)\simeq{\rm
C},\end{equation} that is, for \((f,\lambda),\,(g,\mu)\in
r^{\ast}WZ(S^3)\) and for
\((f',\lambda'),\,(g',\mu')\in (r')^{\ast}WZ((S^3)')\) such that
\(r(f)=r'(f')\) and \(r(g)=r'(g')\), we have:
\begin{equation}
<\,(f,\lambda)\ast(g,\mu),\,(f',\lambda')\ast(g',\mu')\,>
=<\lambda,\lambda'>\ast<\mu,\mu'>,\end{equation}
the right-hand side being the product in \(WZ(\phi)\simeq{\rm C}\).

\subsection{ }

Next we define a connection on \(WZ(S^3)\).   
   They are described as follows.   Let \(b\in \Omega^3_0G\) and
\(U_{f'}\) be a coordinate neighborhood described in 2.5.   
 On \(U_{f'}\)
we put
\begin{equation}
\theta_{U_{f'}}(b)(X)=\frac{i}{48\pi^3}\int_{D'}\,tr(h^{-1}dh)^3\,dX,
\end{equation}
for \(h\in D'b\) and \(X\in D'(Lie\,G)\).   We have
\[\theta_{U_{g'}}=
\theta_{U_{f'}}+{(\chi_{U_{f'},U_{g'}})}^{-1}d\chi_{U_{f'},U_{g'}},\]
where \(\chi_{U_{f'},U_{g'}}\) is the transition function of
\(WZ(S^3)\) :
\[\chi_{U_{f'},U_{g'}}(b)=\chi'(h',k'),\]
for \(h'\in D'b\cap V_{f'}\) and \(k'\in D'b\cap V_{g'}\).   
We have a well defined connection \(\theta\) on \(WZ(S^3)\).   
The curvaure of \(\theta\) becomes 
\begin{equation}
F(X,Y)=-\frac{1}{24\pi^3}\int_{S^3}tr (V^2(XdY-YdX)), 
\qquad V=dff^{-1}\vert S^3.
\end{equation}
The calculation for these formula is the same as in 
[ 6, 10, 11 ].

Similarly we have a connection on \(WZ((S^3)')\) represented by
 a formula 
parallel to ( 2.43 ) but integrated on \(D\).   

On the pullback bundle \(r^{\ast}WZ(S^3)\) there is an induced covariant
derivative:  
\[
(r^{\ast}\bigtriangledown)_{X}s(f)=(\bigtriangledown_{r_{\ast}X}r_{\ast}s)(r(f)),\]
where \(r_{\ast}s\) is the section of \(WZ(S^3)\) defined by
\(r_{\ast}s(b)=
s(f)=[f',c'] \in WZ(S^3)_b\) for a ( and any ) \(f\in Db\).   
\(X\) is a vector field on \(D\), hence \(r_{\ast}X\) is a vector
field on \(S^3\).

Similarly the covariant derivative on \(WZ((S^3)')\) is defined.

   The sections \(WZ(D)\) and \(WZ(D')\) are parallel with respect 
to the respective covariant derivation.   This follows almost from
the definitions by virtue of the infinitesimal form of the 
Polyakov-Wiegman formula:
\begin{equation}
\frac{d}{dt}\vert_{t=0}\,
C_5(f\,e^{tX})=\frac{i}{48\pi^3}\int_{S^4}\,tr(f^{-1}df)^3\,dX,\,\,\mbox{for}\,\,
X\in S^4(Lie\,G),\, f\in S^4G.\end{equation}

\begin{prop}
\begin{eqnarray}
\bigtriangledown WZ(D)&=&0\\
\bigtriangledown WZ(D')&=&0\end{eqnarray}
\end{prop}

\begin{rem}

We could consider in the following construction of the WZW
model those line bundles \(WZ_n(S^3)\) associated to the
\(n\)-th sector of  $\Omega^3G$, but for a fixed \(n\). 
However in the sequel we shall restrict our discussion only to
the contractible component 
$\Omega^3_0G$.\end{rem}

\section{ Construction of WZW actions}
\subsection{ }

Let $\Sigma\in {\cal M}_4$.
Then \(\Sigma\) is a conformally flat manifold with 
boundary 
$\partial\Sigma=\Gamma=\bigcup_{i\in I_+}\Gamma_i
\cup\bigcup_{i\in I_-}\Gamma_i$ with \(\Gamma_i\) a parametrized 
round
\(S^3\) in \(\Sigma\).   

For a $i\in I_-\oplus I_+$, the parametrization
defines the map
\(p_i:S^3\longrightarrow
\Gamma_i\), and the map
\(p_i:\Gamma_iG\longrightarrow
\Omega^3G\), which we denote by the same letter.    Then we have the pull-back bundle
of
$WZ(S^3)$ ( resp. \(WZ((S^3)')\) ) by 
$p_i$.   We define
\begin{eqnarray}
WZ(\Gamma_i)&=& p_i^{\ast}WZ(S^3)\quad\mbox{for
$i\in I_-$,}\nonumber\\
WZ(\Gamma_i)&=& p_i^{\ast}WZ((S^3)')\quad\mbox{for
$i\in I_+$},\end{eqnarray}
then we have respectively
\begin{eqnarray}
WZ(\Gamma'_i)&=& p_i^{\ast}WZ((S^3)')\quad\mbox{for
$i\in I_-$,}\nonumber\\
WZ(\Gamma'_i)&=& p_i^{\ast}WZ(S^3)\quad\mbox{for
$i\in I_+$}.
\end{eqnarray}

The line bundle $WZ(\Gamma)$  
is defined by 
\begin{equation} 
WZ(\Gamma)=\otimes_{i\in
I_-}WZ(\Gamma_i)\otimes\otimes_{i\in I_+}
WZ(\Gamma_i).
\end{equation}

Now let \(\alpha: S^3\longrightarrow S^3\) be the 
restriction on \(S^3\) of a 
conformal diffeomorphism on \(S^4\).   
First we suppose that \(\alpha\) preserves the orientation.  
Then, since the transition function \(\chi\) is invariant
under \(\alpha\), the line bundle \(WZ(S^3)\) is invariant
under \(\alpha\).     If \(\alpha\) reverses the orientation
then \(D\) is mapped to \(D'\) and \(\chi\) is changed to
\(\chi'\).   Then \(\alpha^{\ast}WZ(S^3)\) becomes
\(WZ((S^3)')\).    On the other hand the parametrizations 
\(p_i\) are uniquely defined up to composition with conformal 
diffeomorphisms.    Therefore \(WZ(\Gamma)\) is well defined
for the conformal equivalence class of 
\(\Gamma\in{\cal M}\).

The dual of $WZ(\Gamma)$ is 
\begin{equation} 
WZ(\Gamma^{\prime})= \otimes_{i\in
I_-}WZ(\Gamma_i^{\prime})\otimes\otimes_{i\in I_+}
WZ(\Gamma_i^{\prime}),\end{equation}
and the duality; \(WZ(\Gamma)\times
WZ(\Gamma')\longrightarrow{\rm C}\), is given from ( 2.29 )
by:
\begin{eqnarray*}
&&<\,\otimes_{i\in I_-}[f'_i,c'_i]\otimes
\otimes_{i\in I_+}[g_i,d_i]\,,\,
\otimes_{i\in I_-}[f_i,c_i]\otimes
\otimes_{i\in I_+}[g'_i,d'_i]\,>\\ 
&&=\Pi_{i\in I_-}c_ic^{\prime}_i\cdot \Pi_{i\in I_+}
 d_id^{\prime}_i\cdot 
\exp\{-2\pi i\sum_{i\in I_-}C_5(f_i\vee f^{\prime}_i)
-2\pi i\sum_{i\in I_+}C_5(g_i\vee g^{\prime}_i)\}.
\end{eqnarray*}

The above defined \(WZ(\Gamma)\) satisfies axioms {\bf
A1} and {\bf A2}.  
\medskip

\subsection{ }
 
In the following we shall define step by step the
section
 \(WZ(\Sigma)\) of \(r^{\ast}WZ(\Gamma)\) for
any \(\Sigma\in{\cal M}_4\) with the boundary
\(\partial
\Sigma=\Gamma\) and \(r:\Sigma G\longrightarrow \Gamma G\).

We obtain a compact manifold 
$\Sigma^c\in {\cal M}_4$ without boundary by sewing a copy $D_i$ of
$D$ along 
$\Gamma_i$ for $i\in I_-$ and a copy $D^{\prime}_i$ of 
$D^{\prime}$ for $i\in I_+$;
\[\Sigma^c= (\cup_{i\in I_-}D_i)\,
\cup_{\cup_{I_-}\Gamma_i}\,
\,\Sigma\,
\cup_{\cup_{I_+}\Gamma_i}\,
( \cup_{i\in I_+}D_i^{\prime}).\]

For each boundary component \(\Gamma_i\) of \(\Gamma\) the
parametrization 
\(p_i\) is extended to a parametrization
 \(\tilde p_i:D_i\longrightarrow D\) if \(i\in I_-\), and
\(\tilde p_i:D_i\longrightarrow D'\) if \(i\in I_+\).   The
extension is unique up to composition with conformal
transformations, see 1.1.
   
We put
\begin{eqnarray}
WZ(D_i)&=&(\tilde p_i)^{\ast}WZ(D)\\
WZ(D_i^{\prime})&=&(\tilde p_i)^{\ast}WZ(D').
\end{eqnarray}
For \(i\in I_-\), \(WZ(D_i)\) is a section of the pullback
 bundle of \(WZ(\Gamma_i)\) by the restriction map \(r_i:
D_iG\longrightarrow \Gamma_iG\), and \(WZ(D_i^{\prime})\) is
a section of the pullback bundle of \(WZ(\Gamma_i^{\prime})\)
by the restriction map \(r'_i:D'_iG\longrightarrow
\Gamma_iG\) .    Similarly, for \(i\in I_+\), \(WZ(D'_i)\)
defines a section of the pullback line bundle of
\(WZ(\Gamma_i)\) by \(r'_i\), and \(WZ(D_i)\) is a section of 
 \(r_i^{\ast}WZ(\Gamma_i^{\prime})\) .  
 
{\bf 1}\hspace{0,5cm}
Let \(\Sigma_1\in {\cal M}_4\) and suppose that the
compactified space 
\((\Sigma_1)^c\) is simply connected. that is.
\(\Sigma_1\) is a subset of \(S^4\) deleted several discs
\(D_i;i\in I_-\) and \(D'_i;i\in I_+\) with parametrized
boundaries  
\(\Gamma=\cup_{i\in I_-}\Gamma_i\cup \cup_{i\in 
I_+}\Gamma_i\).   Let
\begin{equation}\Phi_1=\otimes_{i\in
I_+}WZ(D_i)\otimes\otimes_{i\in I_-}
WZ(D_i^{\prime}).\end{equation}
\(\Phi_1\) is a section of the pullback bundle of \(WZ(\Gamma')\) by the 
restriction map 
\[\left(\cup_{i\in I_-}D'_i\cup \cup_{i\in
I_+}D_i\right)G\longrightarrow \left(\cup_{i\in I_-}\Gamma_i\cup
\cup_{i\in I_+}\Gamma_i\right)G.\]
Then \(WZ(\Sigma_1)\) is defined by the duality
relation;
 \begin{equation}
<WZ(\Sigma_1),\Phi_1>=WZ(S^4)=1.\end{equation}
In fact, given \(f\in\Sigma_1 G\), take \(f_i\in
D_iG\), \(i\in I_+\), and \(f_i^{\prime}\in D_i^{\prime} G\),
\(i\in I_-\), in such a way that
\(f\vert\Gamma_i=f_i\vert\Gamma_i\), \(i\in I_+\),
and \(f\vert\Gamma_i=f_i^{\prime}\vert\Gamma_i\), \(i\in I_-\).   
\noindent Let \(WZ(D_i)(f_i)=(f_i,u_i)\), \(i\in I_+\), and
\(WZ(D'_j)(f'_j)=(f'_j,u'_j)\), \(j\in I_-\).   By the definition 
\[u_i\in
WZ(\Gamma'_i)_{r_i(f_i)}, \quad\mbox{and}\quad u'_j\in
WZ(\Gamma'_j)_{r'_j(f'_j)}.\] Then 
\(\Phi_1((f_i)_{i\in I_+},(f'_j)_{j\in
I_-})=((f_i)_{i\in I_+},(f'_j)_{j\in
I_-},\,\otimes_{i\in I_+}u_i\otimes\otimes_{j\in
I_-}u'_j)\).   There is a
\[v\in \otimes_{i\in
I_+}WZ(\Gamma_i)_{r_i(f_i)}\otimes\otimes_{j\in
I_-}WZ(\Gamma_j)_{r'_j(f'_j)}\,=\,WZ(\Gamma)_{r(f)}.\] such that 
\(<v,\,\otimes_{i\in I_+}u_i\otimes\otimes_{j\in
I_-}u'_j>=1. \)  
The definition of \(WZ(D_i)\) and \(WZ(D'_i)\) imply that
\(v\) is independent of  
 \(\{f_i, \,f'_i\}\), but depends only on \(f\).

Thus \(WZ(\Sigma_1)(f)=(f,v)\) is well defined as a section of the
pullback bundle of \(WZ(\Gamma)\) by \(r:\Sigma_1G \longrightarrow
\Gamma G\).

{\bf 2}\hspace{0,5cm}
 Let \(\Sigma_0=S^3\times[0,1]\).   We define 
\begin{equation}
WZ(\Sigma_0)=1_{WZ((S^3)')\otimes WZ(S^3)}.\end{equation}
Then we have
\[ < WZ(\Sigma_0), WZ(D)\otimes
WZ(D')>=<WZ(D),WZ(D')>=1.\]
This is concordant with the definition in
paragraph {\bf 1 }.

{\bf 3}\hspace{0,5cm}
We shall call a \(\Sigma_1\in {\cal M}_4\) described
 in {\bf 1}  that is not of cylinder type a {\it basic
component}.
Any \(\Sigma\in {\cal M}_4\) can be decomposed to a
sum  of several basic components that are patched
together by their parametrized boundaries:
\begin{equation}
\Sigma=\cup_{k=1}^N \Sigma_k.\end{equation}
The incoming boundaries of \(\Sigma_k\) coincide
respectively with the outgoing boundaries of
\(\Sigma_{k-1}\) up to their orientations , that is, \(\Gamma^{k-1}_i=(\Gamma^k_i)'\),
and
\(\Sigma\) is obtained by patching together these
boundaries.   Then there is a duality of
\(WZ(\Gamma^{k-1}_i)= (p^{k-1}_i)^{\ast}WZ((S^3)')\)
and
\(WZ(\Gamma^k_i)= (p^k_i)^{\ast}WZ(S^3)\).   
 Using a suitable Morse
function on \(\Sigma\), we may suppose that the parametrized
boundaries
\(\Gamma_i;\,i\in I_-\) of \(\Sigma\) are all contained
in the boundary \(\partial\Sigma_1\) and
\(\Gamma_i;\,i\in I_+\) are in 
\(\partial\Sigma_N\).   Then we define 
\begin{equation}
WZ(\Sigma_2\cup \Sigma_1)=
<WZ(\Sigma_2),WZ(\Sigma_1)>.\end{equation}
Here \(<\,,\,>\) is the natural pairing 
( contraction ) between the line bundles
\(\otimes_{i\in I_-}WZ(\Gamma_i)
\otimes\otimes_{j\in J_+^1}WZ(\Gamma_j)\) 
and
\(\otimes_{j\in J_-^2}
WZ(\Gamma_j)\otimes\otimes_{k\in J_+^2}
WZ(\Gamma_k)\) .   Here we have written 
\(\partial\Sigma_1=\cup_{j\in
J_+^1}\Gamma_j\bigcup\cup_{i\in
I_-}\Gamma_i\) and \(\partial\Sigma_2=\cup_{k\in
J_+^2}\Gamma_k\bigcup\cup_{j\in
J_-^2}\Gamma_j\), hence 
\(
WZ(\Sigma_2\cup \Sigma_1)\) is a section
of the pullback line bundle of
 \(\otimes_{i\in
I_-}WZ(\Gamma_i)\otimes\otimes_{k\in J_+^2}
WZ(\Gamma_k)\)
 by the boundary restriction map 
 \[
r:\,(\Sigma_2\cup
\Sigma_1)G\longrightarrow 
\cup_{i\in
I_-}\Gamma_iG\cup\cup_{k\in J_+^2}\Gamma_kG,\]
see the explanation after {\bf A3} of 1.2.

\begin{lem}
Let
$\Sigma=\Sigma_1\cup\Sigma_2\cup\Sigma_3$.   
Let their boundaries be 
\[\partial\Sigma_1=\gamma_1\cup\Gamma'_2
\cup\Gamma_3
,
  \quad\partial\Sigma_2=\gamma_2\cup\Gamma_3^{\prime}
\cup\Gamma_1\quad\mbox{and}\quad
\partial\Sigma_3=\gamma_3\cup\Gamma_1^{\prime}
\cup\Gamma_2.\]
Then we have
\begin{equation}<\,<WZ(\Sigma_1),
WZ(\Sigma_2)>,\,WZ(\Sigma_3)\,>=
<\,WZ(\Sigma_1),
<WZ(\Sigma_2),WZ(\Sigma_3)>\,>
\end{equation}
\end{lem}
This is merely the problem of forming a tensor
product of several line bundles, that is a commutative
operation. 

By virtue of this lemma we can form successively 
\begin{eqnarray}
&&WZ(\Sigma_k\cup
\Sigma_{(k-1)}\cdots\cup\Sigma_1)\nonumber\\
&&=
<WZ(\Sigma_k),<WZ(\Sigma_{(k-1)}),
\cdots,WZ(\Sigma_1)>\cdots>.\quad
\end{eqnarray}
This is independent of the order of partition and 
is also independent of how to decompose \(\Sigma^{(k)}=\Sigma_k\cup
\Sigma_{(k-1)}\cdots\cup\Sigma_1\), but depends only on \(\Sigma^{(k)}\).   
Therefore
\[WZ(\Sigma)= WZ(\Sigma_N\cup
\Sigma_{(N-1)}\cdots\cup\Sigma_1)\]
is well defined as a section of the pullback line bundle of 
\(\otimes_{i\in
I_-}WZ(\Gamma_i)\otimes\otimes_{i\in I_+}
WZ(\Gamma_i)\) by the boundary restriction map.

From the construction \(WZ(\Sigma)\) satisfies axiom {\bf A3}.

Now let  \(\Sigma\in{\cal M}_4\) be compact without boundary.   Let \(\Sigma_1\)
and
\(\Sigma_2\) be the basic components such that 
\(\Sigma=\Sigma_1\cup_{\Gamma}\Sigma_2\).  Suppose that 
\[\partial\Sigma_1=\Gamma=\cup_{i\in I}\Gamma'_i\,,\quad
\partial\Sigma_2=\Gamma=\cup_{i\in I}\Gamma_i\,
.
\] Then
from the definition of \(WZ(\Sigma_i)\), \(i=1,2\), we see that:
\begin{eqnarray}WZ(\Sigma)&=&<WZ(\Sigma_2),WZ(\Sigma_1)>=<\otimes_{i\in
I}WZ(D'_i),\otimes_{i\in I}WZ(D_i)>\nonumber \\
\nonumber \\
&=&\sum_{i\in I}1\,.\nonumber\end{eqnarray}

   Thus we have the following

\begin{prop}   For any \(\Sigma\in{\cal M}_4\) which
is compact without boundary 
 \(WZ(\Sigma)\) is a positive integer.
\end{prop}

\begin{prop}
Let \(\Sigma\in{\cal M}_4\) and let \(\Sigma^{ij}\) be
obtained from \(\Sigma\) by identifying the boundaries
\(\Gamma_i\),\(i\in I_-\), and \(\Gamma_j\),\(j\in I_+\),
via
\(\,p_j\cdot(p_i)^{-1}:\Gamma_i\longrightarrow\Gamma_j\).  
Then
\begin{equation}
WZ(\Sigma^{ij})=Tr_{ij}\,WZ(\Sigma),\end{equation}
where \(Tr_{ij}\) are the trace maps ( contraction )
between \(r^{\ast}WZ(\Gamma'_i)\) and 
\(r^{\ast}WZ(\Gamma_j)\) in the tensor product
\( \otimes_{k\in I_-}WZ(\Gamma'_k)\otimes \otimes_{l\in
I_+}WZ(\Gamma_l)\).\end{prop}

\medskip

Connections on \(WZ(\Gamma_i)\) and \(WZ(\Gamma'_i)\) are 
defined
naturally  as the induced connections by ( 3.1 ) and ( 3.2 ).  
Obviously \(WZ(D_i)\) and \(WZ(D'_i)\) are parallel with 
respect to
these connections.   By the formulas of definitions ( 3.3 ), 
( 3.8) and ( 3.13 ) we have a naturally induced connection on
\(WZ(\Gamma)\) with respect to which \(WZ(\Sigma)\) is
 parallel.   
Therefore axiom  {\bf A4} is verified.  

\begin{rem}
 Let $\Sigma\in {\cal M}_4$ and the 
boundary \(\Gamma=\partial\Sigma\) be such that 
$\Gamma=\bigcup_{i\in I_+}\Gamma_i
\cup\bigcup_{i\in I_-}\Gamma'_i$ with \(\Gamma_i\) a 
parametrized 
round \(S^3\).   Let \(r_{\pm}\) denote respectively the restriction maps
onto \(\otimes_{i\in I_{\pm}}(\Gamma_iG)\) .   Then
\begin{equation}
WZ(\Sigma):\,r_-^{\ast}\left(\otimes_{i\in
I_-}WZ(\Gamma_i)\right)\longrightarrow r_+^{\ast}\left(\otimes_{i\in
I_+}WZ(\Gamma_i)\right).\nonumber
\end{equation}
\(WZ(\Sigma)(f)\) for \(f\in \Sigma G\) is the 
higher-dimensional parallel
 transport along the "path" \(f\), [ 19 ].   
When \(I_+=\phi\) or \(I_-=\phi\) we call \(WZ(\Sigma)(f)\) the higher-dimensional holonomy along
\(f\).
\end{rem}

\subsection{ }

To see that the functor \(WZ\) satisfies the axioms of WZW model 
it remains for us to verify the axiom {\bf A5}, the Polyakov-Wiegmann formula on every \(\Sigma\in {\cal M}_4\).    
We have already seen the Polyakov-Wiegmann formula 
 on \(S^4G\), \(DG\) and \(D^{\prime}G\) in ( 2.17 ), ( 2.37 ) and ( 2.40 )
respectively.

Let \(\Sigma\in{\cal M}_4\) with parametrized
boundaries  
\(\Gamma=\cup_{i\in I_-}\Gamma_i\cup \cup_{j\in 
I_+}\Gamma_j\,\).    We shall use the same notation as in 3.1 and 3.2.   
Then the product on each pullback line bundle 
 \(r_i^{\ast}WZ(\Gamma_i)\) ,  \((r^{\prime})_i^{\ast}WZ(\Gamma_i)\), 
 \(r_i^{\ast}WZ(\Gamma^{\prime}_i)\) and  \((r^{\prime})_i^{\ast}WZ(\Gamma^{\prime}_i)\) 
is defined in an obvious manner, and 
the non-vanishing sections  \(WZ(D_i)\) and  \(WZ(D'_i)\)
for \(i\in I_\pm\) 
satisfy the respective Polyakov-Wiegmann formula \[WZ(D_i)(fg)=WZ(D_i)(f)\ast WZ(D_i)(g),\qquad\mbox{etc.}.\] 
  
The products on the line bundle  \(\,S=\otimes_{j\in I_-}(r_j^{\prime})^{\ast}WZ(\Gamma_j)
\otimes \otimes_{i\in I_+}(r_i)^{\ast}WZ(\Gamma_i)\) 
and on the line bundle 
\(\,S^{\ast}=\otimes_{j\in I_-}(r^{\prime}_j)^{\ast}WZ(\Gamma^{\prime}_j)
\otimes \otimes_{i\in I_+}(r_i)^{\ast}WZ(\Gamma^{\prime}_i)\)
are defined by tensoring the product on each \(r_i^{\ast}WZ(\Gamma_i)\), etc..    
We note also that the products are compatible with the duality:
\[<\alpha\ast\beta,\lambda\ast\mu>=<\alpha,\lambda>\ast<\beta,\mu>,\]
for \(\alpha\), \(\beta\in S\) and \(\lambda\), \(\mu\in S^{\ast}\).   Where the product in the right side is that in \(WZ(\phi)\simeq {\rm C}\).  

Now suppose that  \(\Sigma\) is  
a subset of \(S^4\) deleted several discs \(D_i\), \(i\in I_{\pm}\) .   
Let \(r:\,\Sigma G\longrightarrow \Gamma G\) be the restriction map. 
  Then the product on \(r^{\ast}WZ(\Gamma)\) is derived from the product on
 \(S\).   In fact, if we write \(r(f)=(\,r_i(f_i);\,i\in
I_+,\,r^{\prime}_j(f^{\prime}_j);\,j\in I_-\,)\) as in the argument of 3.2, 
then 
\(WZ(\Gamma)_{r(f)}=S_{ r^{\prime}_j(f^{\prime}_j), r_i(f_i)}\),  
so the product on \(S\) yields that on \(r^{\ast}WZ(\Gamma)\), which is seen to
be  independent of the choice of \(\{f_i,\,f^{\prime}_j\}\).     

Let    
\(\Phi_1=\otimes_{i\in I_+}WZ(D_i)\otimes\otimes_{j\in I_-}WZ(D^{\prime}_j)\) .   
\(\Phi_1\) is a section of the line bundle \(S^{\ast}\)  and 
satisfies  
\[\Phi_1(f^{\prime}g^{\prime})=\Phi_1(f^{\prime})\ast\Phi_1(g^{\prime}),\]
for \(f^{\prime}\,,g^{\prime}\in \otimes_{i\in I_+}D_iG\otimes\otimes_{i\in I_-}D^{\prime}_iG\) .   
Since the section \(WZ(\Sigma)\) of \(r^{\ast}WZ(\Gamma)\) was 
defined by the duality;  \(
<WZ(\Sigma),\Phi_1>=WZ(S^4)\),
 we have 
\begin{eqnarray}&&<WZ(\Sigma)(fg),\Phi_1(f^{\prime}g^{\prime})>
=WZ(S^4)(fg\vee f^{\prime}g^{\prime})\nonumber \\ 
&&\quad=WZ(S^4)(f\vee f^{\prime})\ast
WZ(S^4)(g\vee g^{\prime})
=<WZ(\Sigma)(f),\Phi_1(f^{\prime})>\ast
<WZ(\Sigma)(g),\Phi_1(g^{\prime})>\nonumber \\
&&\quad=<WZ(\Sigma)(f)\ast WZ(\Sigma)(g),
\Phi_1(f^{\prime}g^{\prime})>,\nonumber\end{eqnarray}
for any \(f,g\in \Sigma G\) and for \(f^{\prime}\) and  \(g^{\prime}\) that are extensions of \(f\) and \(g\) to \(\cup_{i\in
I_-}D'_i\cup\cup_{j\in I_+}D_j\) respectively.  
Therefore we have 
\[WZ(\Sigma)(fg)=WZ(\Sigma)(f)\ast WZ(\Sigma)(g),\]
for \(f,\,g\in \Sigma G\).

 Let
\(\Sigma=\Sigma_1\cup_{\Gamma}\Sigma_2\).   The
product operations on \((r_i)^{\ast}WZ(\Gamma_i)\), \(i=1,2\), are
compatible with the contraction, in particular we have
\begin{eqnarray}
<WZ(\Sigma_1)(f_1)\ast WZ(\Sigma_1)(g_1)&,&\,
WZ(\Sigma_2)(f_2)\ast
WZ(\Sigma_2)(g_2)>\nonumber\\
&=&WZ(\Sigma)(f)\ast WZ(\Sigma)(g),
\end{eqnarray}
where \(f,g\in \Sigma G\) and \(f_i=f\vert\Sigma_1\),
\(i=1,2\) etc..    For a general
\(\Sigma\in{\cal M}_4\) the formula follows from  ( 3.14 ) and the
definition of
\(WZ(\Sigma)\) in ( 3.12 ).   Thus we have proven the following generalization of the Polyakov-Wiegmann formula.

\begin{thm}
\begin{equation}
WZ(\Sigma)(f)\ast WZ(\Sigma)(g)=WZ(\Sigma)(fg)
\end{equation}
for \(f,g\in\Sigma G\).
\end{thm}

\section{ Extensions of the group \(\Omega^3_0G\)}

It is a well known observation
that the two-dimensional WZW action gives a geometric
description of the central extension \(\widehat{LG}\) of the loop group \(L G\).   
The associated group cocycle yields a Lie algebra cocycle for the affine Kac-Moody
algebra based on \(Lie (G)\),  [ 2, 7 ].     The total space of the
\(U(1)\)-principal bundle  \(\widehat{LG}\) was described as the set of
equivalence classes of pairs \((f,c)\in D^2G\times U(1)\),   where \(D^2\) is the
2-dimensional disc with boundary \(S^1\).    The equivalence relation was defined
on the basis of Polyakov-Wiegmann formula [ 13 ], as it was so in our
four-dimensional generalization treated in section 2.     

\noindent Associated to the line
bundle \(WZ(S^3)\) there exists a \(U(1)\)-principal bundle over 
\(\Omega^3_0G\).   However this bundle has not any natural group structure contrary 
to the case of the extension of loop group.      Instead 
J. Mickelsson in [ 11 ] gave an extension of \(\Omega^3_0G\) by 
the abelian group \(Map({\cal A}_3,U(1))\), where \({\cal A}_3\) is the space of connections on \(S^3\) .       
In the following we shall explain after [ 10 ] two extensions of 
  \(\Omega^3_0G\) by 
the abelian group \(Map({\cal A}_3,U(1))\) that are in duality.

\subsection{ \,}

We consider the quotient space 
\begin{equation}\widehat{\Omega G}=D^{\prime}G\times Map({\cal A}_3,U(1))/\sim^{\prime},
\end{equation}
where \(\sim^{\prime}\) is the equivalence relation defined by
\begin{equation}
(f^{\prime},\lambda)\sim^{\prime} (g^{\prime},\mu)
 \mbox{ if and only if }
\left\{
\begin{array}{rl} f'\vert S^3&=g'\vert S^3\\
\mu(A) &=\lambda(A)\chi'(f',g')\quad \mbox{ for any} A\in {\cal A}_3.
\end{array}
\right. \end{equation}

The projection \(\pi: \widehat{\Omega G}\longrightarrow \Omega^3_0G\) is defined by \(\pi([f^{\prime},\lambda])= f^{\prime}|S^3\).   Then \(\widehat{\Omega G}\) becomes a principal bundle over \(\Omega^3_0G\) with the structure group \(Map({\cal A}_3,U(1))\).   Here the \(U(1)\) valued transition function \(\chi'(f',g')\) is considered as a constant function in \(Map({\cal A}_3,U(1))\).

The group structure of  \(\widehat{\Omega G}\) is given by the MickelssonÕs 
2-cocycle ( 2.9 ) on \(D^{\prime}\):
\[
\gamma_{D^{\prime}}(\cdot;f^{\prime},g^{\prime}),
\qquad \mbox{for} \,f^{\prime},g^{\prime}\in D^{\prime}G.\]
We note that, since it is the coboundary of 
\[\frac{i}{24\pi^3}\int_{D^{\prime}}\alpha_4(A;f^{\prime}),\]
\(\gamma_{D^{\prime}}\) is in fact a cocycle.   
We define the product on \(D^{\prime}G\times Map({\cal A}_3,U(1))\) by
\begin{equation}
(f^{\prime},\lambda)\ast(g^{\prime},\mu)=\left(f^{\prime}g^{\prime}
,\lambda(\cdot)\mu_{f^{\prime}}(\cdot)\exp\{2\pi i\gamma_{D^{\prime}}(A;f^{\prime},g^{\prime})\}\,\right),
\end{equation}
where \[
\mu_{f^{\prime}}(A)=\mu((f^{\prime}|S^3)^{-1}A(f^{\prime}|S^3)+(f^{\prime}|S^3)^{-1}d(f^{\prime}|S^3)).\]
Then \(D^{\prime}G\times Map({\cal A}_3,U(1))\) is endowed with a group structure and \(\widehat{\Omega G}\) inherits it.   The group \(Map({\cal A}_3,U(1))\) is embedded as a normal subgroup in   \(\widehat{\Omega G}\).   Thus  \(\widehat{\Omega G}\)  is  an extension of \(\Omega^3_0G\) by the abelian group  \(Map({\cal A}_3,U(1))\) ,  [ 10, 11 ].      

We have another extension of  \(\Omega^3_0G\) by  \(Map({\cal A}_3,U(1))\) if
 we consider 
\begin{equation}
\widehat{\Omega^{\prime} G}=DG\times Map({\cal A}_3,U(1))/\sim,
\end{equation}
where the equivalence relation \(\sim\) is defined by
\begin{equation}
(f,\lambda)\sim (g,\mu)
 \mbox{ if and only if }
\left\{
\begin{array}{rl} f\vert S^3&=g\vert S^3\\
\mu(A) &=\lambda(A)\chi(f,g)\quad \mbox{ for any} A\in {\cal A}_3.
\end{array}
\right. \end{equation}
The product on \(\widehat{\Omega^{\prime} G}\) is defined by the
 same way as above using the 2-cocycle 
\(\gamma_D(A;f,g)\) of ( 2.9 ), and \(\widehat{\Omega^{\prime} G}\) becomes a
extension of \(\Omega^3_0G\) by the abelian group 
 \(Map({\cal A}_3,U(1))\).   

The group \(Map({\cal A}_3, U(1))\) acts on \({\rm C}\) by 
\(\lambda\cdot c=\lambda(0)c\).   Then the associated line bundle to \(\widehat{\Omega G}\) is \(WZ(S^3)\) and that associated to 
 \(\widehat{\Omega^{\prime} G}\) is \(WZ((S^3)^{\prime})\).

\begin{rem}
Consider the empty three manifold \(\phi\) and look it as the  boundary of \(S^4\) .   Then we  may follow the above definition to have an extension of \(\phi G\) by 
\(Map({\cal A}_3, U(1))\).   It becomes  \(\widehat{\phi G}=S^4G\times Map({\cal A}_3, U(1))/ \sim\), where 
\[(F,\lambda)\sim(G,\mu)\mbox{ if and only if }
 \mu(A) =\exp\{2\pi i\omega(F,F^{-1}G)\}\lambda(A)\quad\mbox{for
any \(A\).}\] Then, since \((F,\lambda)\sim(F,\lambda(0))\), it
reduces to 
\(\widehat{\phi G}=S^4G\times  U(1) /\sim\),  that is, \(\widehat{\phi G}\simeq U(1)\).   
The product in \(\widehat{\phi G}\) may be defined by the same formula as in
 ( 4.3 ), but we have seen that it reduces to that of ( 2.16 ) because of the
equality 
\(\gamma_{S^4}(A;F,G)=\gamma(F,G)\), ( 2.10 ).
\end{rem}

The duality between two extensions \(\widehat{\Omega G}\)
and \(\widehat{\Omega^{\prime} G}\) is given as follows.
For \([f^{\prime},\lambda]\in \widehat{\Omega G}\) and  \([f,\alpha]\in \widehat{\Omega^{\prime} G}\),  we put
\begin{equation}
<[f^{\prime},\lambda], [f,\alpha]>=[f\vee f ^{\prime},
\lambda(0)\alpha(0)],
\end{equation}
where on the right hand side we used the product in \(\widehat{\phi G}\simeq U(1)\).
In fact, suppose that  \((f^{\prime},\lambda)\sim^{\prime}(g^{\prime},\mu)\) and  \((f,\alpha)\sim(g,\beta)\) .   Then we have
\begin{eqnarray}
\mu(A)\beta(A)&=&\lambda(A)\alpha(A)\exp\{2\pi i\left[
\gamma_{D^{\prime}}(A;f^{\prime},g^{\prime})+\gamma_D(A;f,g)\right]\}\\
&=&\lambda(A)\alpha(A)\exp\{2\pi i
\gamma(f\vee f^{\prime},g\vee g^{\prime})\}.
\end{eqnarray}
Here we used the relation ( 2.11 ).

The Lie algebra cocycle corresponding to the group cocycle
\(\gamma_{D}\) is calculated in [ 11 ].    It is given by
\begin{eqnarray}
c(A;X,Y)&= &\frac{i}{12\pi^2}\int_{D}tr\,dA(dXdY+dYdX)\nonumber\\
&=&\frac{i}{12\pi^2}\int_{S^3}\,tr(A(dXdY+dYdX).
\end{eqnarray}

The Lie algebra cocycle corresponding to the group cocycle 
\(\gamma_{D^{\prime}}\) is given by
\begin{eqnarray}
&&\frac{i}{12\pi^2}\int_{D^{\prime}}tr\,dA(dXdY+dYdX)\nonumber\\
&&=-\frac{i}{12\pi^2}\int_{S^3}\,tr(A(dXdY+dYdX)=-c(A;X,Y)
\end{eqnarray}

\subsection{ Remarks}
\noindent{\bf 1}
\hspace{0.5cm}
The Euclidean action of a field
\(\varphi: \Sigma\longrightarrow G\) in WZW confomal field
theory is defined as 
\begin{equation}
S_{\Sigma}(\varphi)=-\frac{ik}{12\pi^2}\int_{\Sigma}
tr(d\varphi^{-1}\wedge\ast d\varphi)+C_{\Sigma}(\varphi).\end{equation}
 \(S_{\Sigma}(\varphi)\)
is invariant under a conformal change of metric and 
the second term \(C_{\Sigma}(\varphi)\) is required to obtain a
conformal invariance of the action.   This was shown by K.
Fujii in [ 6 ], and first noticed by E. Witten in [ 20 ] for the 
two-dimensional WZW model.   The kinetic term in ( 4.11 ) is linear
with respect to the multiplication of the fields;
\begin{equation}
\int_{\Sigma}
tr(d(fg)^{-1}\wedge\ast d(fg))=\int_{\Sigma}
tr(df^{-1}\wedge\ast df)+\int_{\Sigma}
tr(dg^{-1}\wedge\ast dg)
\end{equation}
and does not affect the Polyakov-Wiegmann formula.   
Hence we prefered only to deal with the topological term
\(C_{\Sigma}(f)\), [ 3 ].

\noindent{\bf 2}\hspace{0.5cm}
  The argument in this paper will be valid also for \(2n\)-dimensional conformally 
flat manifolds with boundary if the Lie group \(G=SU(N)\) is such that \(N\geq n+1\), in
this case we have \(\pi_{2n}(G)=0\) and \(\pi_{2n+1}(G)={\bf Z}\).     We shall have also
the abelian extensions of \(\Omega^{2n-1}_0G\) by
\(Map({\cal A}_{2n-1},U(1))\).    For that purpose we must have the Polyakov-Wiegmann
formula for the action functional
  \[C_{2n+1}(f)=\frac{-i}{(2n-1)!(2\pi i)^{2(n-1)}}
  \int_{D^{2n+1}}\,tr(\tilde g^{-1}d\tilde g)^{2n+1},\quad g\in S^{2n}G .\]
  See [ 6 ].     It seems that Polyakov-Wiegmann formula has not
yet been proved for general \(n\) larger than 3.
  
\noindent{\bf 3} 
\hspace{0.5cm} Losev, Moore, Nekrasov and Shatashvili
[ 9 ] discussed a four-dimensional WZW theory based 
on K\"ahler manifolds.   Their Lagrangian is defined
by 
\[-\frac{1}{4\pi}\int_{\Sigma}\omega\wedge
Tr(g^{-1}\partial g\wedge \ast
g^{-1}\overline{\partial}g)+
\frac{i}{12\pi}\int_{\Sigma \times [0,1]}\omega\wedge
Tr(g^{-1}dg)^3.\]
The theory has the finiteness properties for the
one-loop renormalization of the vacuum state.   The
authors studied the algebraic sector of their 
theory.   The category of algebraic manifolds is
not well behaved under contraction, hence their
theory does not fit our axiomatic description.    

\noindent{\bf 4}
\hspace{0.5cm} \(S^4\) is obtained by patching together
two quarternion spaces and we have the conjugation
\(q\longrightarrow q^{-1}\) on it.   Under the
conjugation 
\(WZ(S^4)\) is invariant but \(WZ(D)\) and \(WZ(D')\) will
interchange.   Since the conjugation inverts the orientation, \(WZ(\Sigma)\) is
invariant under the conjugation of
\(\Sigma\).   We can convince ourselves of this
fact if we follow the argument to define
\(WZ(\Sigma)\) for a \(\Sigma\in 
{\cal M}_4\) .    This is the CPT
invariance.

\medskip
{\bf Acknowledgements}

I would like to express my thanks to Professor J. Mickelsson of Royal Institute of Technology in
Sweden for his valuable comments and criticisms which helped me several times to avoid incorrect arguments.     
  I also thank Dr. Y.Terashima of Tokyo University for his remarks about the 
precise formulation of the functor \(WZ\) . 
I wish to acknowledge Professor M. Guest of Tokyo
Metropolitan University who gave valuable comments on an earlier version.

\end{document}